\documentclass[a4, 12pt]{amsart}
\usepackage{graphicx,epsfig}
\usepackage{amssymb,amsfonts,amsmath,amsthm,dsfont,wasysym,pifont,stmaryrd,  xy}
\input xy
\xyoption{all}
\usepackage{epstopdf,yfonts}
\usepackage{hyperref}
\usepackage{pdflscape}
\usepackage[vcentermath]{youngtab}
\numberwithin{equation}{section}

\newtheorem{thm}{Theorem}[section]
\newtheorem{prop}[thm]{Proposition}

\newtheorem{lem}[thm]{Lemma}

\theoremstyle{definition}

\newtheorem*{remark*}{Remark}

\newcommand{\cO}{{\mathcal O}}

\newcommand{\cS}{{\mathcal{S}}}

\newcommand{\tr}{\text{Tr}}
\newcommand{\fp}{\kappa}

\newcommand{\lm}{\mathfrak{L}}

\newcommand{\C}{\mathbf{C}}

\newcommand{\ind}{\mathrm{Ind}}

\newcommand{\glO}[2]{\text{GL}_{#1}({\cO}_{#2})}
\newcommand{\slO}[2]{\text{SL}_{#1}({\cO}_{#2})}
\newcommand{\sslO}[2]{\text{sl}_{#1}({\cO}_{#2})}
\newcommand{\spO}[2]{\text{Sp}_{#1}({\cO}_{#2})}

\newcommand{\uO}[2]{\text{U}_{#1}({\cO}_{#2})}
\renewcommand{\O}[2]{\text{O}_{#1}({\cO}_{#2})}

\newcommand{\mO}[2]{\text{M}_{#1}({\cO}_{#2})}
\newcommand{\irr}{\mathrm{Irr}}
\newcommand{\res}{\mathrm{Res}}
\begin{document}
\title[Representations of Classical Groups]{On representations of classical groups over principal ideal local rings of length two}
\author{Pooja Singla}
\address{ Center for Advanced Studies in Mathematics \\
Department of mathematics \\
Ben Gurion University of the Negev\\
P.O.B. 653 \\
 Beer Sheva 84105\\
 Israel 
}
\email{pooja@math.bgu.ac.il}
\begin{abstract} We study the complex irreducible representations of special linear, symplectic, orthogonal and unitary groups over principal ideal local rings of length two. We construct a canonical correspondence between the irreducible representations of all such groups that preserves dimensions.  
The case for general linear groups has already been proved by author. 
\end{abstract}
\maketitle
\section{Introduction} 
Let $\mathrm{F}$ be a non-Archimedean local field with ring of integers $\cO$.
Let $\wp$ be the unique maximal ideal of $\cO$ and $\pi$ be a fixed
uniformizer of $\wp$.
Assume that the residue field $\cO/\wp$ has odd characteristic $p$.
We denote by $\cO_{\ell}$
the reduction of $\cO$ modulo $\wp^{\ell}$, i.e., $\cO_{\ell} =
\cO/\wp^{\ell}$. Therefore $\cO_1$ will denote the residue field of $\cO$. 

  The representation theory of classical groups over the rings $\cO$ and the finite rings $\cO_{\ell}$ has attracted attention of many mathematcians. See Singla~\cite{MR2684153} for the history of this problem for the General and Special linear groups over the rings $\cO$ and their finite quotients $\cO_{\ell}$. In the direction of the classical groups, Lusztig~\cite{MR2048585} constructed several irreducible representations of reductive groups over finite rings and Jaikin-Zapirain~\cite{MR2169043} looked at the problem of constructing irreducible representations of compact pro-p groups. But
the knowledge of all irreducible representations of classical groups over $\cO$ and the finite rings $\cO_{\ell}$ is still far from complete.         
 
 In Singla~\cite{MR2684153}, we gave a method to construct irreducible representations of general linear groups $\glO{n}{2}$. In this article, we use it to construct irreducible representations of other classical groups over the rings $\cO_2$. The questions of this article are also motivated by a conjecture of Onn~\cite[Conjecture 1.2]{MR2456275}, 
which says that the isomorphism type of the group algebra of automorphism group of finite $\cO$-modules depend on the ring of integers only through the cardinality of its residue field. More 
generally one can ask this question for the other classical groups over $\cO_{\ell}$ as well. In Singla~\cite{MR2684153}, we proved Onn's conjecture for groups $\glO{n}{2}$.  
 
\subsection{Main Results} The classical groups over ring $\cO_2$ are defined as: 
\begin{enumerate}

\item Special Linear Group: $\slO{n}{2} = \{ g \in \glO{n}{2} \mid \det(g) = 1 \}$.  \\

\item Symplectic Group:  $\spO{n}{2} = \{ g \in \glO{2n}{2} \mid g^{\mathrm{t}} J g = J \}$, where $g^{\mathrm{t}}$ denotes transpose of $g$ and $J = \tiny {\left[ \begin{matrix} 0 & I_n \\ -I_n & 0 \end{matrix} \right]}. $ \\

\item  Orthogonal Group: $\O{n}{2} = \{ g \in \glO{n}{2} \mid g^{\mathrm{t}} g = I_n \}$. \\

\item Unitary Group: Let $\tilde{F}$ be a degree two unramifield extension of $F$ and $\sigma$ be the unique nontrivial Galois automorphism of $\tilde{F}$. 
Let 
$\tilde{\cO}_{2}$ be the corresponding unramified extension of $\cO_2$, then $\sigma$ restricts to an automorphism of $\tilde{\cO}_2$ (denoted again as $\sigma$).  
Applying $\sigma$ entry wise we obtain an automorphism of $\mathrm{GL}_n(\tilde{\cO}_2)$. For any $g \in \mathrm{GL}_n(\tilde{\cO_2})$, let $g^{\star} = (g^{\sigma})^{\mathrm{t}}$. 
Then the unitary group is defined as   
$$\uO{n}{2} = \{g \in \mathrm{GL}_n(\tilde{\cO}_2) \mid g g^{\star}	 = I_n \}. $$
\end{enumerate}
We prove that the number and dimensions of irreducible representations of these classical groups over $\cO_2$ depend on $\cO$ only through the cardinality of residue field. 
More precisely we prove the following.

Let $\mathrm{F}$ and $\mathrm{F}'$ be local fields with rings of integers
  $\mathcal{O}$ and $\mathcal{O'}$, respectively, such that their residue fields
  are finite and isomorphic (with a fixed isomorphism). Let $\wp$ and $\wp'$ be the
  maximal ideals of $\mathcal{O}$ and $\mathcal{O}'$ respectively. As
  described earlier, $\cO_2$ and $\cO'_2$ denote the rings
  $\cO/\wp^2$ and $\cO'/\wp'^2$, respectively. For a ring $R$, we use
$C(R)$ as collective notation for any of the classical groups $\mathrm{SL}_n(R)$, $\mathrm{Sp}_n(R)$, $\mathrm{O}_n(R)$, or $\mathrm{U}_n(R)$ over $R$.   

\begin{thm}
\label{thm:main} 
There exists a canonical bijection between the irreducible representations of
  $C(\cO_2)$ and those of 
  $C(\cO'_2)$, which preserves dimensions. 

\end{thm}
By the equivalence between the number of conjugacy classes and distinct irreducible representations, we also obtain that the number of conjugacy classes of 
the classical groups over $\cO_2$ depend on $\cO$ only through $|\cO_1|$. We remark that very little is known about the conjugacy classes of classical group $C(\cO_{\ell})$ (See~\cite{MJUP, MR2543507, mythesis}).

\section{Proof of Theorem~\ref{thm:main}} First of all, we set up few notations.
By character we shall always mean one-dimensional representation. For an abelian group $A$, the set of its characters is denoted by $\hat{A}$. For any group $G$, the set of inequivalent irreducible representations is denoted by $\mathrm{Irr}G$. 
Let $$\kappa : \glO{n}{2} \rightarrow \glO{n}{1} \,\,
\mathrm{and}\,\, \bar{\kappa} : C(\cO_2) \rightarrow C(\cO_1)$$ 
be the natural quotient maps with $K = \ker(\kappa)$ and
$L(C) = \ker(\bar{\kappa})$. 
We shall use the following results of Clifford theory. 

\begin{thm}[Clifford Theory]
\label{thm:clifford} Let $G$ be a finite group and $N$ be a normal subgroup. Let $\rho$ be an 
  irreducible representation of $N$ and $T(\rho) =
  \{ g \in G \mid \rho^{g} = \rho \}$ be the stabilizer of $\rho$. Then
  the following hold 
\begin{enumerate}

\item If $\pi$ is an irreducible representation of $G$ such that $\langle
  \pi|_{N}, \rho \rangle \neq 0$, then $\pi|_{N} = e (
  \oplus _{\rho \in \Omega} \rho )$ where $\Omega$ is an orbit of irreducible
  representations of $N$ under the action of $G$, and $e$
  is a positive integer. \\

\item Let $A = \{ \theta \in \irr (T(\rho)) \mid \langle \res ^{T(\rho)}_{N} \theta, \rho \rangle
\neq 0 \}$ and $$ B = \{ \pi \in \irr G \mid \langle \res^{G}_{N} \pi, \rho \rangle
\neq 0 \}.$$
Then $$ \theta \rightarrow \ind^{G}_{T(\rho)}(\theta)$$ is a bijection
of $A$ onto $B$. \\

\item Let $H$ be a subgroup of $G$ containing $N$, and suppose that $\rho$ has an extension $\tilde{\rho}$ to
$H$ (i.e., $\tilde{\rho}|_{N} = \rho $). Then the representations $\chi \otimes \tilde{\rho}$
for $\chi \in \irr(H/N)$ are irreducible, distinct for distinct $\chi$, and \\
\[
\ind^{H}_{N}(\rho) = \oplus_{\chi \in \irr(H/N)} \chi \otimes
\tilde{\rho}.
\]   
\end{enumerate}    
\end{thm}
Let $\widehat{L(C)}$ denote the set of characters of $L(C)$. The group $C(\cO_2)$ acts on $\widehat{L(C)}$ by conjugation. That is, if $\alpha \in C(\cO_2)$ and $\phi \in \widehat{L(C)}$ then $\phi^{\alpha}(x) = \phi(\alpha x \alpha^{-1})$ for $x \in L(C)$.
For any $\phi \in \widehat{L(C)}$, let $T_C(\phi) = \{ \alpha \in C(\cO_2) \mid \phi^{\alpha} = \phi \}$ be the stabilizer of $\phi$ in $C(\cO_2)$. 
\begin{prop}
\label{prop: sl2} For any  $\phi \in \widehat{L(C)}$, 
\begin{enumerate} 
\item[(a)] There exists a canonical character $\chi_{\phi}$ of $T_C(\phi)$ such that 
$\chi_{\phi}|_{L(C)} = \phi$. 
\item[(b)] The group $T_C(\phi)/K$ depends on $\cO$ only through $|\cO_1|$.
\item[(c)] The cardinality $|C(\cO_2)/T_C(\phi)|$ depends on $\cO$ only through $|\cO_1|$. 
\end{enumerate}
\end{prop} 
We postpone the proof of this proposition to \S~\ref{sec: class} and ~\ref{sec: sl(p/n)}. Assuming this, we complete the proof of Theorem~\ref{thm:main}.

\begin{proof}[Proof of Theorem~\ref{thm:main}] Let $\mathcal{S}_C$ denote 
the set of $C(\cO_2)$-orbits in 
$\widehat{L(C)}$ under conjugation. Therefore, by Clifford Theory, there exists a bijection (also canonical by proposition~\ref{prop: sl2}) between the sets
\begin{eqnarray}
\label{correspondence}
 \coprod_{\phi \in \mathcal{S}_C} \{ \mathrm{Irr}(T_C(\phi)/K) \} &
\longleftrightarrow & \rm{Irr}C(\cO_2),
\end{eqnarray}
given by, 
\begin{eqnarray}
\delta \mapsto \mathrm{Ind}_{T_C(\phi)}^{C(\cO_2)}(\chi_{\phi} \otimes \tilde{\delta}),
\end{eqnarray}
where $\tilde{\delta}$ is a representation of $T_{C}(\phi)$ obtained by composing $\delta$ with the natural projection map $T_C(\phi) \rightarrow T_C(\phi)/K$. Since the left side of (\ref{correspondence}) depends on $\cO$ only through the cadinality of $\cO_1$, theorem follows. 

\end{proof}
\section{The groups $\O{n}{2}$, $\spO{n}{2}$, $\uO{n}{2}$, and $\slO{n}{2}(p \nmid n)$}
\label{sec: class}
Fix a nontrivial additive character $\psi: \cO_1 \rightarrow \mathbf C^{\times}$.
For each 
$A \in M_n(\cO_1)$ the character $\psi_A : K \rightarrow \mathbf C^{\times}$ is defined by 
\[
\psi_A(I + \pi X) = \psi(\tr (AX)).
\]
The bilinear form $\langle.,.\rangle : M_n(\cO_1) \times M_n(\cO_1) \rightarrow \cO_1$, defined by $(A,B) \mapsto \tr(AB)$ is non-degenerate therefore the assignment $A \mapsto \psi_A$ defines an isomorphism $M_{n}(\cO_1) \cong \hat{K}$. Let $T_G(\psi_A) = \{ g \in \glO{n}{2} \mid \psi_A^g = \psi_A \}$ (in place of $T(\psi_A)$ of Singla~\cite{MR2684153} to remove any ambiguity) denote the stabilizer of $\psi_A$ in
$\glO{n}{2}$. In this section we prove Proposition~\ref{prop: sl2} for Orthogonal, Unitary, Symplectic and Special linear groups$(p\nmid n)$. For this section $C(\cO_2)$ denotes either $\O{n}{2}$, $\uO{n}{2}$, $\spO{n}{2}$ or $\slO{n}{2} (p \nmid n)$. For any classical group $C(\cO_2)$, let $M_C$ be the subgroup
of $M_n(\cO_1)$ such that $X \mapsto I + \pi X$ defines an isomorphism $M_C \cong L(C)$. By restricting $\langle .,. \rangle$ to $M_C$, we obtain a bilinear form on $M_C$ as well. Observe that if this restriction is non-degenerate then,    
\begin{enumerate}
\item The map $A \mapsto \psi_A|_{L(C)}$ defines an isomorphism, $M_c \cong \widehat{L(C)}$.
\item $ T_C(\psi_A|_{L(C)}) = T_G(\psi_A) \cap \C(\cO_2)$. 
\end{enumerate}
  
We recall the following result of Singla~\cite{MR2684153}  
\begin{prop}
\label{prop:main} For any $\phi \in \hat{K}$, there exists a 
canonical character $\chi_{\phi}^G$ of $T_G(\phi)$ such that $\chi_{\phi}^G|_{K} = \phi$ (such a 
character $\chi_{\phi}^G$ is called an extension of $\phi$). 
\end{prop}
\begin{proof} 
For proof see Proposition~2.2 and Section~5 of Singla~\cite{MR2684153}.
\end{proof}

 Let $\phi = \psi_A|_{L(C)}$ for some $A \in M_C$, then part (a) of Proposition~\ref{prop: sl2} follows by taking $\chi_{\phi} = \chi_{\phi}^G|_{T_C(\phi)}$. Parts (b) and (c) of Proposition~\ref{prop: sl2} follow by the following facts:
\begin{enumerate}
\item  $T_C(\phi)/L(C) \cong Z_{\glO{n}{1}}(A) \cap C(\cO_1)$. 
\item $|C(\cO_2)/T_C(\phi)| = |C(\cO_1)/ (Z_{\glO{n}{1}}(A) \cap C(\cO_1))|$.   
\end{enumerate}
These follow easily by definitions of $T_C(\phi)$, $L(C)$ and Corollary~5.3 of Singla~\cite{MR2684153}.
  
We shall show that the bilinear form $\langle .,. \rangle|_{M_C}$ is non-degenerate for Symplectic, Unitary, Orthogonal, and Special linear($p \nmid n$) groups. Further in \S~\ref{sec: sl(p/n)} we show that it is not true for Special linear group with $p \mid n$, which makes this case bit harder. We shall use a completely different method to solve $\slO{n}{2}(p \mid n)$.   

\subsection{$\O{n}{2}$} The kernel $L(O)$ of the natural projection map $\O{n}{2} \rightarrow \O{n}{1}$ is isomorphic to  
$$M_O = \{ X \in \mO{n}{1} \mid X + X^{\mathrm{t}} = 0 \}$$ by $I + \pi X \mapsto X$. 
For any $A \in M_O$, 
since $\langle.,.\rangle$ is non-degenerate on $\mO{n}{1}$, there exists $ Y \in \mO{n}{1}$ such that $\tr(AY) \neq 0$. Take $B = Y - Y^{\mathrm{t}} \in M_{O}$. By using additivity and $\tr(Z) = \tr(Z^{\mathrm{t}})$ for any $Z \in \mO{n}{1}$, we obtain 
\[
\tr(AB) = \tr(A(Y-Y^{\mathrm{t}})) = \tr (AY + A^{\mathrm{t}} Y^{\mathrm{t}}) = 2 \tr(AY) \neq 0. 
\]
Here we have also used the fact that the residue field is of odd characteristic. Hence $\langle.,.\rangle|_{M_O}$ is non-degenerate.
\subsection{$\uO{n}{2}$} The kernel $L(U)$ of the natural projection map $\uO{n}{2} \rightarrow \uO{n}{1}$ is isomorphic to the set 
$$M_U = \{ X \in \mO{n}{1} \mid X + X^{\star} = 0 \}$$ by $I + \pi X \mapsto X$. The rest of the argument for this follows similar to orthogonal group case, by replacing $Y - Y^{\mathrm{t}}$ with $ Y - Y^{\star}$.  
\subsection{$\spO{n}{2}$}
\label{sym}
The kernel $L(\mathrm{Sp})$ of the natural projection map $\spO{n}{2} \rightarrow \spO{n}{1}$ consists of matrices $I + \pi X$, $X \in M_{2n}(\cO_2)$ such that 
\[
(I + \pi X)^{\mathrm{t}} J (I + \pi X) = J, \; \mathrm{where}\; J = \left[ \begin{matrix} 0 & I_n \\ -I_n & 0 \end{matrix} \right]. 
\]
Let $$M_{\mathrm{Sp}} = \{ X \in \mO{2n}{1} \mid X^{\mathrm{t}} J + J X = 0 \} $$ 
then $X \mapsto I + \pi X$ is easily seen to give an isomorphism between $M_{\mathrm{Sp}}$ and $L(\mathrm{Sp})$. 
Also observe that, $$M_{\mathrm{Sp}} = \{ \begin{bmatrix} U & V \\  W & -U^{\mathrm{t}} \end{bmatrix} \mid U, V, W \in \mO{n}{1}, V = - V^{\mathrm{t}}, W = W^{\mathrm{t}} \}.$$ 
For any $A = \begin{bmatrix} U & V \\  W & -U^{\mathrm{t}} \end{bmatrix}, \, X = \begin{bmatrix} U' & V' \\ W' & -U'^{\mathrm{t}} \end{bmatrix} \in M_{\mathrm{Sp}}$, 
$$\tr(AX) = 2 \tr(UU') + \tr(VW') + \tr(WV').$$ By using the fact that $(A, B) \mapsto \tr(AB)$ is a non-degenerate bilinear form on the set of matrices, symmetric and skew-symmetric 
matrices over $\cO_1$. We obtain that for any nonzero $A \in M_{\mathrm{Sp}}$ there exists 
$X \in  M_{\mathrm{Sp}}$ such that $\tr(AX) \neq 0$. Therefore $\langle.,.\rangle|_{M_{\mathrm{Sp}}}$ is non-degenerate.

\subsection{$\slO{n}{2}$, $p \nmid n$:}
\label{sec: sl(p-n)} Let $L(SL)$ denote the kernel of the natural projection map $\slO{n}{2} \rightarrow \slO{n}{1}$. Identify the set $L(\mathrm{SL})$ with 
$M_{\mathrm{SL}} = \{ X \in \mO{n}{1} \mid \tr(X) = 0 \}$, by $I + \pi X \mapsto X$. We prove that the bilinear form $(A,B) \rightarrow \tr(AB)$ is non-degenerate 
on $M_{\mathrm{SL}}$. 

For any non-diagonal matrix $A = (a_{ij}) \in \sslO{n}{1}$, there exists pair $(i_0,j_0)$ such that
$i_0 \neq j_0$ and $a_{i_0 j_0} \neq 0$. Take $ B = (b_{ij}) \in M_{\mathrm{SL}}$ such that $b_{j_0 i_0} = 1$ and zeros everywhere else. Then $\tr(AB) \neq 0$. On the other hand if
$A = (a_{ij}) \in M_{\mathrm{SL}}$ is a non-zero diagonal matrix, then $p \nmid n$ implies there exists $i_0 \neq j_0$ such that $a_{i_0i_0} \neq a_{j_0j_0}$. Then by taking $B = (b_{ij}) \in M_{\mathrm{SL}}$ to be the matrix satisfying 
$b_{i_0i_0} = 1$, $b_{j_0j_0} = -1$ and $b_{ij} = 0 $ for all $(i,j) \notin\{(i_0, i_0), (j_0, j_0)\}$, we obtain that $\tr(AB) \neq 0$. This proves the assertion.   
 
\section{Special Linear Group $\slO{n}{2}$, $p \mid n$}
\label{sec: sl(p/n)} Let $L(SL)$ denote the kernel of the natural projection map $\slO{n}{2} \rightarrow \slO{n}{1}$. As mentioned in \S~\ref{sec: sl(p-n)}, the set $L(\mathrm{SL})$ can be identified with 
$M_{\mathrm{SL}} = \{ X \in \mO{n}{1} \mid \tr(X) = 0 \}$, by $I + \pi X \mapsto X$. We firstly show that $\langle.,.\rangle|_{M_{\mathrm{SL}}}$ is not non-degenerate by showing that the scalar matrices lie in its radical. Let $A = aI_n \in \mO{n}{1}$, then $p \mid n$ implies  
$A \in M_{\mathrm{SL}}$. Hence for any $X \in M_{\mathrm{SL}}$ we obtain $\tr(AX) = a\tr(X) = 0$. This implies that $\langle.,.\rangle|_{M_{\mathrm{SL}}}$ is not non-degenerate. Hence the method discussed in the last section does not work in this case. \\  

Define an equivalence relation on $M_n(\cO_1)$ by $A \sim B $ if there exists a scalar
$x \in \cO_1$ such that $A = xI + B$. Denote the equivalence class of $A$ under this relation by $[A]$ and the set of 
equivalence classes by $\lm$. The set $\lm$ is an abelian group under the operation $[A] + [B] = [A+B]$.  
For any $[A] \in \lm$ define $\psi_{[A]} : L(SL) \rightarrow \mathbf C^{\times}$ by
\[ 
\psi_{[A]} (I + \pi X) = \psi (\text{Tr}(AX)).     
\]
Then $\psi_{[A]}$ is well defined character of $L(SL)$ and $[A] \mapsto \psi_{[A]}$ gives an isomorphism 
$\lm \tilde{\rightarrow} \widehat{L(SL)}$. The group $\glO{n}{2}$ acts on $\lm$ by conjugation
via its quotient $\glO{n}{1}$, and therefore on $\widehat{L(SL)}$. For $\alpha \in \glO{n}{2}$ and 
$\psi_{[A]} \in \widehat{L(SL)}$,
we obtain,  
\[
\begin{matrix}
\psi^{\alpha}_{[A]}(I + \pi X) & =  &  \psi(\tr(A \fp(\alpha) X \fp(\alpha)^{-1})) &  = &  \psi_{\fp(\alpha)^{-1}[A] \fp(\alpha)} (I + \pi X). 
\end{matrix}
\]
Let $T_G(\psi_{[A]}) = \{ \alpha \in \glO{n}{2} \mid \psi_{[A]}^{\alpha} = \psi_{[A]} \}$. 
Observe that $L(SL)$ is a subgroup of $I + \pi M_n(\cO_2)$ and the character $\psi_{A} \in \widehat{K}$ restricts to 
$\psi_{[A]}$ on $L(SL)$. By definitions it follows that $T_G(\psi_A) = \{ \alpha \in \glO{n}{2} \mid \psi_A^\alpha = \psi_A \}$ is a subgroup of $T_G(\psi_{[A]})$. Let $T_{\mathrm{SL}}(\psi_{[A]})$ be the stabilizer of $\psi_{[A]}$ 
in $\slO{n}{2}$, then $T_{\mathrm{SL}}(\psi_{[A]}) = T_{G}(\psi_{[A]}) \cap \slO{n}{2}$. 
We subdivide our further discussion to two cases.

\subsection{ The case $T_G(\psi_A) = T_G(\psi_{[A]})$:} The condition $T_G(\psi_A) = T_G(\psi_{[A]})$ implies $$T_{\mathrm{SL}}(\psi_{[A]}) = T_G(\psi_{A}) \cap \slO{n}{2}.$$
Then define $\chi_{\psi_{[A]}} = \chi_{\psi_A}^G|_{T_G(\psi_{A}) \cap \slO{n}{2}}$, where $ \chi_{\psi_A}^G$ is as obtained from Proposition~\ref{prop:main}. Then $\chi_{\psi_{[A]}}|_{L(SL)} = \psi_{[A]}$. This proves the existence of canonical extension, that is Proposition~\ref{prop: sl2}(a), in this case.  

\subsection{The case $T_G(\psi_{[A]}) \neq T_G(\psi_{A})$:} Let $T_{\mathrm{SL}}(\psi_A)$ denote the set $T_G(\psi_A) \cap \slO{n}{2}$. For this case, we follow the following steps to obtain the character $\chi_{\psi_{[A]}}$ of $T_{\mathrm{SL}}(\psi_{[A]})$.\\

\noindent {\bf Step~1:} We show that the group $T_{\mathrm{SL}}(\psi_A)$ is normal in $T_{\mathrm{SL}}(\psi_{[A]})$ and the group  $T_{\mathrm{SL}}(\psi_{[A]})/T_{\mathrm{SL}}(\psi_A)$ is abelian.\\ 
{\bf Step~2:} There exists an abelian group $\mathcal{S}$ such that 
\[ T_{\mathrm{SL}}(\psi_{[A]}) = T_{\mathrm{SL}}(\psi_A)\mathcal{S},
\]
and the intersection $\mathcal{S} \cap T_{\mathrm{SL}}(\psi_A) $ is trivial. \\
{\bf Step~3} There exists a character ${\chi}_A$ of $T_{\mathrm{SL}}(\psi_A)$ that is invariant in $T_{\mathrm{SL}}(\psi_{[A]})$ and hence extends to give required character $\chi_{\psi_{[A]}}$. \\ 

By definition of the action of $\glO{n}{2}$ on $\mO{n}{1}$ and the set $\{[X] \mid X \in \mO{n}{1}\}$ of equivalence classes, we obtain $T_G(\psi_{[A]}) = \{ g \in \glO{n}{2} \mid g[A]g^{-1} = [A]\}$ and $T_G(\psi_A) = \{ g \in \glO{n}{2} \mid g A g^{-1} = A \}$. Let $g \in T_G(\psi_{[A]})$ be such that $gAg^{-1} = A + xI$ for some $x \in \cO_1$, and let $z \in T(\psi_A)$. Then, 
\[
\begin{matrix} (gzg^{-1}) A (gzg^{-1})^{-1} & = & gzg^{-1} A gz^{-1}g^{-1} & = & gz(A-xI)z^{-1} g^{-1} = A.          
\end{matrix} 
\]
This implies that $T_G(\psi_A) (T_{\mathrm{SL}}(\psi_A))$ is a normal subgroup of $T_G(\psi_{[A]})(T_{\mathrm{SL}}(\psi_{[A]}))$. Further the quotient $T_{G}(\psi_{[A]})/T_G(\psi_A)(T_{\mathrm{SL}}(\psi_{[A]})/T_{\mathrm{SL}}(\psi_A))$ is abelian because \[g_1 g_2 A (g_1 g_2)^{-1} = (g_2 g_1) A (g_2 g_1)^{-1}. \] This completes the proof of Step~1.  

We can assume that matrix $A$ has all its eigenvalues in the field $\cO_1$, for if not we can apply the argument to an extension field of $F$. Hence for further discussion, we shall assume that that all the eigenvalues of $A$ lie in the field $\cO_1$.
 
The assumption $T_G(\psi_{[A]}) \neq T_G(\psi_A)$ implies there exists a nonzero scalar $x \in \cO_1$ such that $A$ is conjugate to $xI +A$. Therefore, if $a$ is an eigenvalue of $A$ then so is $a+x$. We arrange the distinct eigenvalues of $A$ in the following order,  
\[
a_{11}, a_{12}, \ldots, a_{1p}, a_{21}, a_{22}, \ldots, a_{2p}, \ldots, a_{r1}, a_{r2}, \ldots, a_{rp},
\]  
where $a_{ij} = a_{i(j-1)} + x$, $a_{ip} + x = a_{i1}$ for all $1 \leq i \leq p$,$ 2 \leq j \leq p$, and for $i \neq i'$, $a_{ij}-a_{i'j} \notin (x)$ (the additive 
space generated by $x$). Assume that $A$ is in its Jordan Canonical form (see Theorem~3.5 of Singla~\cite{MR2684153}), that is 
\begin{equation}
\label{jordan}
A = \oplus_{i=1}^r \oplus_{j=1}^p A_{ij},
\end{equation} 
where each $A_{ij}$ is further direct sum of Jordan blocks with unique eigenvalue $a_{ij}$. 

Every element of $T_G(\psi_{[A]})$ modulo the group $T_G(\psi_A)$ just permutes the matrices $A_{ij}$ among each other. Therefore every element of $T_G(\psi_{[A]})$ can be written as product of an element of the permutation matrix and an element of the group $T_G(\psi_A)$. The cosets of $T_{\mathrm{SL}}(\psi_A)$ in $T_{\mathrm{SL}}(\psi_{[A]})$ can be parametrized by the permutation matrices consisting of $pr \times pr$ blocks 
with each $(i,j)^{th}$ block of size equal to size of matrix $A_{ij}$ and with the property that each block is either identity or zero matrix.
Let $\mathcal{S}$ be the collection of the permutation matrices corresponding to the quotient $T_{\mathrm{SL}}(\psi_{[A]})/T_{\mathrm{SL}}(\psi_A)$, then $\cS$ is an abelian group (by Step~1)) satisfying, 
\begin{enumerate}
\item $T_{\mathrm{SL}} = \cS T_{\mathrm{SL}}(\psi_A).$ 
\item The intersection $\cS \cap T_{\mathrm{SL}}(\psi_A)$ is trivial. 
\end{enumerate}
This proves Step~2. For step~3, first of all we briefly recall the construction of character $\chi^G_{\psi_A}$ of $T_G(\psi_A)$ that extends $\psi_A$, from Singla~\cite{MR2684153}. \\

Let $s: \cO_1^{\times} \rightarrow \cO_2^{\times}$ be
the unique multiplicative section of the natural projection map $\cO_2^{\times} \rightarrow \cO_1^{\times}$. By defining $s(0) = 0$, we obtain a section $s: \cO_1 \rightarrow \cO_2$ of the natural projection map $\cO_2 \rightarrow \cO_1$. By 
extending it entry wise, we obtain a map from $\glO{n}{1} \rightarrow \glO{n}{2}$, denoted again by $s$. 

For a matrix $A$ as given in (\ref{jordan}), let $Z_{\glO{n}{2}}(s(A)) = \{ g \in \glO{n}{2} \mid gs(A) = s(A)g \}$ and $m_{ij}$ be the size of each matrix $A_{ij}$. Then $$Z_{\glO{n}{2}}(s(A)) = \Pi_{i=1}^r \Pi_{j=1}^p Z_{\glO{m_{ij}}{2}}(s(A_{ij})),$$ and $T_G(\psi_A) = K.Z_{\glO{n}{2}}(s(A))$ (see Lemma~5.1 and Corollary~5.3 of Singla~\cite{MR2684153}). 

For any $a \in \cO_1$ define a character $\psi_a: 1+ \pi \cO_2 \rightarrow \mathbf C^{\times}$ by
$\psi_a(1+ \pi x) = \psi(ax)$. Since $\cO_2^{\times}$ is direct product of $\cO_1^{\times}$ and $1+\pi \cO_2$, the characters 
$\psi_a$ extend trivially to $\cO_2^{\times}$, denote this by $\chi_a$. Define a character $\chi$ of 
$Z_{\glO{n}{2}}(s(A)) = \Pi_{i=1}^r \Pi_{j=1}^p Z_{\glO{m_{ij}}{2}}(s(A_{ij}))$ by 
$$\chi(\Pi_{i=1}^r \Pi_{j=1}^p X_{ij}) = \chi_{a_{11}}(\det(X_{11}))\ldots \chi_{a_{rp}}(\det(X_{rp})).$$ 
Then the character $\chi_{\psi_A}^G = \psi_A.\chi$ of $T_G(\psi_A)$, defined by $\psi_A.\chi(uv) = \psi_A(u)\chi(v)$ for all $u \in K$ and $v \in Z_{\glO{n}{2}}(s(A))$, satisfies $ \chi_{\psi_A}^G|_{K} = \psi_A$. 
Let ${\chi}_A = \chi_{\psi_A}^G|_{T_{\mathrm{SL}}}(\psi_A)$. Then,  
\begin{lem} The one dimensional representation ${\chi}_{A}$ of $T_{\mathrm{SL}}(\psi_A) $ is fixed by $T_{\mathrm{SL}}(\psi_{[A]})$.
\end{lem}
\begin{proof} To prove that ${\chi}_A$ is fixed by $T_{\mathrm{SL}}(\psi_A)$, it is 
sufficient to prove that the restriction of $\chi$ to $Z_{\glO{n}{2}}(s(A)) \cap \slO{n}{2}$ is invariant under the action of elements of the group $\cS$, as $\psi_A|_{K \cap \slO{n}{2}} = \psi_{[A]}$ and $\psi_{[A]}$ is fixed by $T_{\mathrm{SL}}(\psi_{[A]})$ by definition of $T_{\mathrm{SL}}(\psi_{[A]})$.

Observe that any permutation matrix $s \in \mathcal{S}$ such that $sAs^{-1} = A + x'I$ permutes the matrices $A_{ij}$ and $A_{i'j'}$ among 
each other only if both $A_{ij}$ and $A_{i'j'}$ have the same block decomposition and $a_{ij} - a_{i'j'} \in (x')$, where $(x')$ denotes the additive 
space generated by $x'$. Let $X \in (Z_{\glO{n}{2}}(s(A)) \cap \slO{n}{2})$. Then by the structure of $Z_{\glO{n}{2}}(s(A))$, we obtain that
$X = \Pi_{i=1}^r \Pi_{j=1}^p X_{ij}$. Let $\det(X_{ij})= \beta_{ij}(1+\pi \alpha_{ij})$, where $\beta_{ij} \in \cO_1^{\times}$ and $1+\pi \alpha_{ij} \in 1 + \pi \cO_2$. Here we have again used the fact that $\cO_2^{\times}$ is direct product of $\cO_1^{\times}$ and $1 + \pi \cO_2$. Then $ X \in \slO{n}{2}$ implies that $\Pi_{i=1}^r \Pi_{j=1}^p \beta_{ij} = 1$ and therefore  
$\sum_{i=1}^r \sum_{j=1}^p \alpha_{ij} = 0$. By definition of $\chi$, we obtain 
$$\chi(X) = \psi(\sum_{i=1}^r \sum_{j=1}^p (a_{ij} + (j-1)x)(\alpha_{ij}))$$ 
and $$\chi(sXs^{-1}) = \psi( \sum_{i=1}^r \sum_{j=1}^p (a_{ij} + x' + (j-1)x)(\alpha_{ij})).$$ but then
$\sum_{i=1}^r \sum_{j=1}^p \alpha_{ij} = 0$ implies 
$$\sum_{i=1}^r \sum_{j=1}^p (a_{ij} + (j-1)x)(\alpha_{ij}) = \sum_{i=1}^r \sum_{j=1}^p (a_{ij} + x' + (j-1)x)(\alpha_{ij}). $$
Therefore $\chi(X) = \chi(sXs^{-1})$. This proves the lemma. 
\end{proof}

Hence by Singla~\cite[Lemma~5.4]{MR2684153} (also see Isaacs~\cite[Exercise~6.18]{MR0460423}) 
$\chi_A$ extends to $T_{\mathrm{SL}}(\psi_{[A]})$. Furthermore, as $\cS$ is an abelian group with a trivial intersection with $T_{\mathrm{SL}}(\psi_A)$, we can define extension $\chi_{\psi_{[A]}}$ canonically such that $\chi_{\psi_{[A]}}|_{\cS}$ is trivial. This completes the proof of Proposition~\ref{prop: sl2}(a) for $\slO{n}{2}(p \mid n)$.\\

Let $Z_{\glO{n}{1}}([A])$ be the subgroup of $\glO{n}{1}$ consisting of matrices $g \in \glO{n}{1}$ satisfying $g[A]g^{-1} = [A]$. Parts (b) and (c) of Proposition~\ref{prop: sl2} follow by observing,
\begin{enumerate} 
\item The group $T_{\mathrm{SL}}(\psi_{[A]})/L(SL)$ is isomorphic to 
$Z_{\glO{n}{1}}([A]) \cap\slO{n}{1} $.
\item $|\slO{n}{2}/ T_{\mathrm{SL}}(\psi_{[A]})| = |\slO{n}{1}/(Z_{\glO{n}{1}}([A]) \cap \slO{n}{1})|$. 
\end{enumerate}

\subsection{Acknowledgments} The author is greatly thankful to Uri Onn for suggesting these questions and for many useful comments on a preliminary version of this article. It is also pleasure to thank Amritanshu Prasad for 
very useful feedback. The author was partially supported by Center for Advanced Studies in Mathematics at Ben Gurion University, Beer Sheva, Israel.

\end{document}